\documentclass{notices}
\usepackage{amsfonts,amssymb,amsmath,amscd}
\usepackage{graphicx}
\usepackage{subcaption,float}
\usepackage{xurl}
\usepackage{tabularx}
\usepackage{wrapfig}
\usepackage{makecell}

\title{Are University Budget Cuts Becoming A Threat to Mathematics? with Additional Discussion}

\author{Edgar Fuller
  \affil{Edgar Fuller is Director of the Center for Transforming Teaching in Mathematics and Distinguished University Professor of Mathematics at Florida International University. He served as Chair of the Department of Mathematics at West Virginia University from 2008-2018. He can be reached at efuller@fiu.edu.}
}

\begin{document}
\maketitle

\section*{Introduction}
Mathematics as an area of study occupies an important place in higher education. Due in part to its utility in other disciplines as well as its role in student learning, institutions of higher education (IHEs) often have large numbers of mathematics faculty with different balances of teaching and research in different ranks and appointment structures. Most flagship IHEs, especially state land-grant institutions, have large undergraduate populations taking mathematics courses in many cases built around the widespread use of calculus~\cite{bressoud2012second} and the connections between mathematics and science, technology, and engineering. These connections have made mathematics departments essential to universities~\cite{olson2012engage} and emphasized the critical role math plays in supporting student success~\cites{reinholz2020time,calcscience} in all areas of post-secondary education. We tend to take that essential nature of mathematics at the undergraduate level, and for research universities at the graduate level, as a given, but that characterization no longer holds for some IHEs.

In September of $2023$, the Board of Governors for West Virginia University (WVU) voted to discontinue the graduate program in mathematics, ending a period of $35$ years of doctoral studies in mathematics and a much longer history of providing Masters level opportunities to students at the institution. The vote resulted in the loss of 16 of the 48 faculty lines in the School of Mathematical and Data Sciences~\cite{wvu2023} including the termination of tenure-track as well as teaching faculty. The changes at WVU, characterized as ``Academic Transformation," have received significant national attention summarized here~\cite{ams2023sm} in more detail. The president and provost at WVU based these cuts on recommendations from the consulting firm rpkGroup~\cite{ukrpk2022} hired in response to declining enrollment. They echo similar recommendations made by the same consulting group to the University of Kansas in $2022$, and represent growing concerns in higher education about projected future decreases in enrollment known as the ``enrollment cliff"~\cite{enrcliff2023}. In Kansas, $42$ academic programs were recommended for discontinuation, with reasoning, as in prior evaluations in Arkansas and Missouri~\cite{ams2023sm}, that questioned the funding of mathematics and other STEM degree programs. In yet another review, the Kansas Board of Regents suggested that mathematics requirements should be removed altogether~\cite{ams2023sm} for some majors. 

\begin{figure}[h!]
\centering
\includegraphics[width=.95\columnwidth]{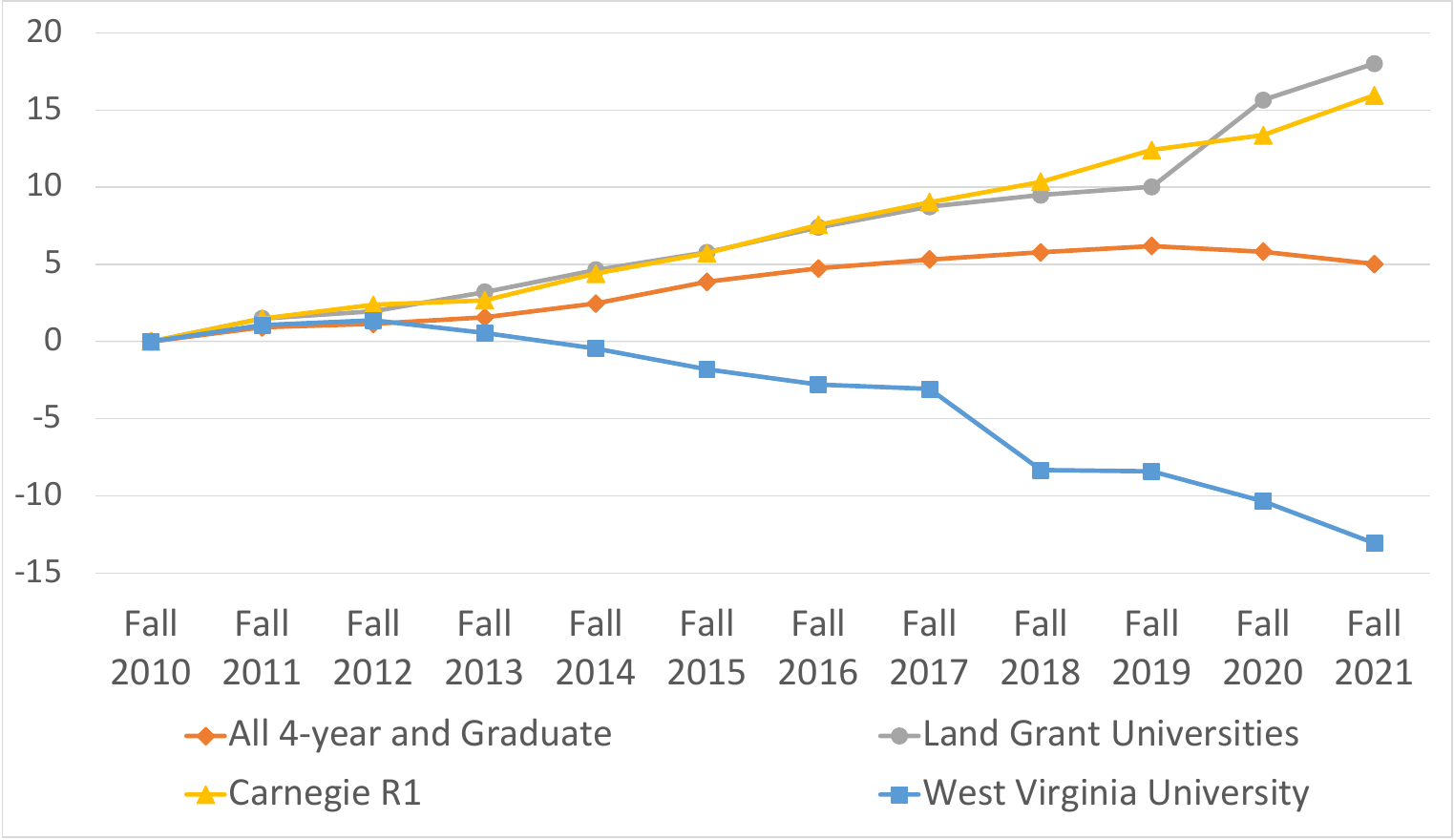}
\caption{Percent change in enrollment since 2010}\label{fig:enroll}
\end{figure}

Data from the Integrated Postsecondary Education Data System (IPEDS)~\cite{ams2023sm} shows that since the overall post-secondary enrollment maximum in $2010$, student populations at both land-grant and research universities have increased as shown in Figure \ref{fig:enroll}, even after experiencing small declines during the Covid-19 lockdown. Total enrollments at four-year institutions have decreased recently, but the decreases at WVU follow neither of these patterns.

While these trends may change, overall student interest in attending large land-grant and research focused institutions remains consistent with the enrollments that supported their programs over the last decade. That broader pattern nationally suggests that factors other than enrollment trends, such as already existing plans to restructure programs, may be responsible for these choices. Even in the context of the decreases in student enrollment at WVU, the removal of the graduate programs in mathematics represents a fundamental change in the mission and values of an institution of higher education. Calls to shift resources away from mathematics education and research in particular ignore the powerful role it plays within an institution and in society.

At WVU, the restructuring runs even deeper. A total of 28 of 338 degree programs, including both the MS and PhD in mathematics~\cite{wvurec} have been discontinued and more than $146$ of the tenured and non-tenured faculty at WVU, or roughly $16\%$ of total faculty university-wide, will either resign by the end of the 2023-24 academic year, retire early during the 2023-24 academic year, or be terminated. Within this group, $44$ tenured or tenure-track faculty will be involuntarily separated at the end of the academic year, leading the American Association of University Professors to state in a recent letter~\cite{ams2023sm} that with these measures enacted ``...tenure, and by extension academic freedom, cannot be said to exist at West Virginia University." At least $112$ of the positions will be removed from academic colleges, or more than $15$\% of the $880$ academic faculty at WVU. Concerns presented by consulting firms~\cite{ams2023sm} focus on the belief that universities should only allocate resources to disciplines with some specific value~\cite{ukrpk2022}, but how that value is measured is often inconsistent and not reflective of the true state of an institution or its mission. Administrators at WVU ignored the fact that tuition revenue generated from mathematics courses averaged \$19.3 million per year over the last three years, making it one of the most efficient departments on campus with operational costs annually around $35$\% of generated revenue. Instead they chose to focus solely on the fact that the total value of grants obtained by faculty was less than $1$ million per year even though the total was above that cutoff as recently as $2018$. Afterwards, President Gordon Gee stated that ``...I don’t believe every aspect of math is essential"~\cite{ams2023sm} and others in the Office of the Provost have also disputed any need for theoretical mathematics. The university has since begun removing mathematics coursework from degree programs in business and other areas.

When universities, or the consultants they hire, overlook the significant components of what mathematics provides to higher education, it signals a fundamental change in the strategic goals of those institutions. The true impact of a discipline like mathematics becomes lost in the narratives used to justify decision making, and these trends begin to affect mathematics departments within large flagship research universities as well as universities which have already seen dramatic resource reallocations. Leaders in mathematics both inside and aligned with higher education must work together to correct these narratives and reinforce the transformative potential of mathematics.

Land grant institutions, for example, were created by the Morrill Act "in order to promote the \textit{liberal and practical education} of the industrial classes in the several pursuits and professions in life."~\cite{morrill1862}. Policymakers at land-grant IHEs may discuss the practical uses of education as the primary driver of programmatic choices, but the use of the word liberal here implies a specific mission as well. The government sought to increase access to education in a way that would broaden the opportunities available to citizens, and increase their freedom to choose new occupational pathways. The act provided a mechanism to create institutions of higher education so that communities would have access to as many areas of learning as possible. Land-grant IHEs in particular represent essential paths to higher education for students in their respective states, a role that should be maintained. Many land-grant and Hispanic Serving Institutions (HSIs) have developed mathematics departments which excel~\cite{garcia2019defining} in research while meeting the needs of the communities they serve providing indispensable access to general education mathematics courses as well as research and graduate study beyond those needs. To overlook the added value of graduate study in mathematics at an IHE to its overall research and education missions, especially at these universities, neglects the critical place they occupy in our communities.  We can't predict when the next Katherine Johnson, a graduate student in mathematics at WVU in 1940~\cite{malcom2020katherine}, will arrive from Appalachia asking to study the subjects that will be used to launch a rocket to the stars, or when the next Addison Fischer, who earned Bachelors and Masters degrees in Mathematics at WVU from 1966-1972, will want to study mathematics and then go on to build multi-billion dollar corporations. IHEs are founded on the notion of providing that access, and the strong support of mathematics amplifies academic quality and breadth across every program. Universities, and the mathematics departments within them, should certainly operate efficiently and provide measurable returns, but support for graduate programs fulfills the higher purpose of making strong educational opportunities broadly accessible to the people we serve. As a profession we must continue to demonstrate this impact, even when consultants lack the expertise to recognize it. We have an obligation to work to make sure that the systems in which we operate see the power of mathematics to elevate and inspire consistently and indisputably. Not long before his passing, Bob Moses summarized the power of mathematics~\cite{bobmoses2021}, saying that
\begin{quote}
\it Math literacy will be a liberation tool for people trying to get out of poverty and the best hope for people trying not to get left behind.
\end{quote}

As a profession, we often view the importance of our discipline as self-evident since the role of mathematics in other disciplines, or just as a practical tool in our lives, generally remains well-established. Reinforcing that importance requires ongoing work, however, since changing priorities, enrollment trends and funding stresses at IHEs create challenges that recur. In the mid 1990s, the University of Rochester proposed the elimination of its mathematics graduate programs due to financial difficulties~\cite{rochester1996} as a part of its ``Renaissance Plan", but the national response from multiple science and higher education organizations led to a reversal of the decision. Years later, the department was thriving~\cite{jackson1997whatever}, giving some indication of departmental strengths that the initial plan overlooked and providing some sense that targeting graduate programs as a cost saving measure may not be the best strategy. The events at WVU and other IHEs highlight a potentially changing landscape in which even core disciplines such as mathematics may be targeted for reduction when not aligned with strategic goals that have departed from the historical mission of universities and colleges. Multiple factors have led higher education to this place, and it is worth considering research and graduate study in mathematics in the context of some of the drivers of that change.

\section*{Changes in Enrollment and State Appropriations}
\begin{figure}[h!]
\centering
\includegraphics[width=.95\columnwidth]{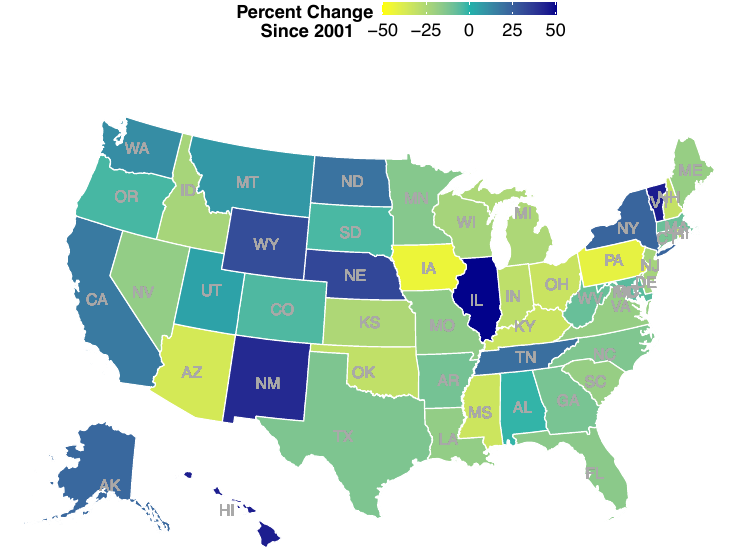}

\caption{Percent change in state appropriated funding for higher education from SHEF report~\cite{sheff2023} by state from 2001 to 2022 in CPI adjusted constant dollars}\label{fig:stateappchg}
\end{figure}

\begin{figure}[h!]
\centering
\includegraphics[width=.95\columnwidth]{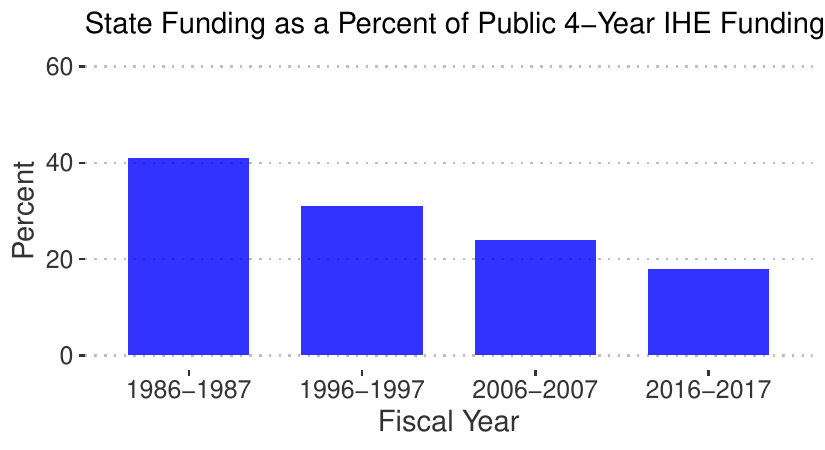}
\caption{Percent of total university budgets for public four year IHEs from state appropriations - APLU}\label{fig:pctbudget}
\end{figure}

Enrollment declines at some IHEs have influenced many of the current discussions of program efficiency, but public IHEs have already been experiencing resource strain in recent years~\cite{sheff2023} as decreasing state appropriations have forced more dependence on tuition and federal funding~\cite{ams2023sm}. Figure \ref{fig:stateappchg} shows the change in state appropriated funding for the operation of public IHEs by state from $2001$ to $2022$. Over the last twenty years, the majority of states have decreased appropriations substantially, in constant dollars over the last twenty years. Data from the Association of Public and Land-grant Universities (APLU), shown in Figure \ref{fig:pctbudget}, illustrates the decline in the state funding allocations to $4$-year public IHEs over the last four decades, where it has decreased to approximately half the level of the late 1980s.

State appropriated funds, particularly for land grant universities, remains an important part of the funding model for public universities~\cite{ams2023sm} though this source has decreased for IHEs in 36 of 50 states for the last two decades, leading universities to increase tuition and turn to other sources of revenue to make up for shortfalls. In WV, for example, state funding has decreased by 8.5\% in that period. As a result, tuition has steadily risen at state institutions, increasing the student share of the cost of attending an IHE~\cite{ams2023sm}.

As the total of state funding has decreased, the appropriated funds can be used in different ways by IHEs as well. From 2001 to 2021 the fraction of state appropriated funding for higher education used for IHE operations nationally remained relatively constant, averaging $78.2\%$, while in WV the fraction specifically used for operations decreased from $65.5\%$ to $56.2\%$. The reallocated $10\%$ of the overall already decreasing state appropriations moved to the research, agricultural extension programs, public health care services, and medical education (RAM) category indicating a use outside of academic operations. In $2022$, WV allocated $18.08\%$ of total higher education appropriations specifically to public health care services and medical education, followed by state funding used to support agriculture extension programs ($5.27\%$), university affiliated hospital support ($1.67\%$), and state support for research ($0.38\%$). Revenue and expense data from WVU~\cite{wvurec} show that medical education programs operate at a loss, and so this shift of state appropriations towards that portion of the university appears to be a strategic effort to support medical activity at the expense of core academic operations. The combined decrease in appropriations overall with allocation to other uses applies stress to the ability of state IHEs to meet mission goals. IHEs then begin to look at reductions to even high enrollment departments such as mathematics that normally are seen as essential to the proper functioning of a university in an effort to balance budgets. The resulting reviews often focus on the cost of a department, major enrollments, or other metrics, but they ignore many of the important impacts of academic departments such as mathematics at an IHE.

\section*{University and Mathematics Research Activity}
As resources have tightened at even larger universities, arguments made in support of mathematics research and graduate programs sometimes focus on the classification of an IHE as a Research - Very High Activity institution, commonly referred to as R1 status. A total of $146$ IHEs are currently classified as R1, or about $3.7\%$, making this an important distinction for some universities. Those within the R1 group form a fairly broad collection of private and public institutions and for those IHEs R1 status is considered a measure of prestige that enhances the reputation of an IHE as well as that of its faculty. Both IHEs and mathematics departments within them appeal at times to this high level research status for recruitment and to procure additional research contracts. These contracts in turn provide additional revenue to universities that can be critically important during times of financial stress, leading some administrators to emphasize grant funded research over other departmental contributions, but R1 status does not necessarily capture the specifics of research at an IHE. 

Grouping IHEs into private, public, Hispanic-serving (HSI) and land-grant (LG) institutions, some differences in the characteristics of the institutions appear in Table \ref{tab:R1char}. Land-grant institutions ($35$\%), which include some public, private, and HSI IHEs, form a large subset of R1s. HSIs which are R1s (11 of 146, $7.5$\%) are a growing subset that includes four land-grant IHEs. A Forbes analysis estimated that HSIs have been underfunded in $16$ states by almost $\$13$ billion in the last 30 years~\cite{forbes2023}, with the mean number of faculty~\cite{ams2023sm} reported as $63.5\%$ of the average number for other R1s as of $2021$ though their undergraduate enrollment is considerably higher. They award a substantial number of PhDs and significantly more undergraduate degrees than other R1s. 

The table summarizes the mean values of several key aspects of R1 IHEs including total undergraduate enrollment (UG), total graduate enrollment (GR), number of Bachelors (BS) and PhDs (PhD) awarded, and the number of full-time faculty in ladder rank (assistant, associate, and full professors). The mean number of students enrolled at the undergraduate level provides the clearest discrepancy between the public and private R1s. External funding listed here as Science and Engineering Research and Development Expenditures  (SERD, in millions of dollars) are considerably higher for private and land-grants, while the average number of faculty and the number of PhDs awarded are almost the same. This metric along with the number of postdoctoral and non-faculty research staff with doctorates (PD) is the largest factor impacting R1 ranking as of the $2021$ analysis. PD numbers are calculated differently at some IHEs, with some counts including faculty equivalent appointments that involve teaching and others large numbers of laboratory support staff with no academic duties.

\begin{table}[h!]
\centering
\footnotesize
\begin{tabular}{p{.2\columnwidth}cccc}
& \thead{\small Private} & \thead{\small Public} & \thead{\small HSI} & \thead{\small LG} \\
 Number       & 39     & 107      & 11       & 51      \\
 SERD   & 579.7  & 388.8   & 235.6   & 483.7   \\
 UG& 8966 & 25241 & 28382 & 26331 \\
 GR& 9315 & 7787  & 6838  & 7792  \\
 Mean PhD     & 331  & 345   & 287   & 415   \\
 Mean BS      & 2629 & 6228  & 6882  & 6700  \\
 Mean FAC     & 1471 & 1474  & 962   & 1624  \\
 Mean PD      & 887  & 443  & 299   & 580   \\
\end{tabular}\caption{Characteristics of R1 universities by group} \label{tab:R1char}
\end{table}

\begin{figure}[h!]
\centering
\includegraphics[width=.9\columnwidth]{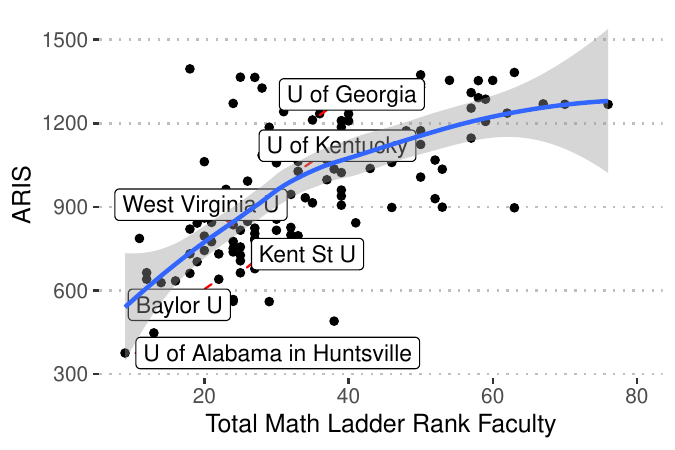}
\caption{Scatterplot of the ARIS weighted ranking by the total of tenure line math faculty}\label{fig:matharisTT}
\end{figure}

As a part of their ``Academic Transformation" process, administrators at WVU identified departments with graduate programs with less than $\$1,000,000$ per year in grant related research expenditures as measured by the SERD data in Table \ref{tab:R1char}, and reviewed student enrollments in their programs. The subsequent evaluation of the School of Mathematical and Data Sciences led to the discontinuation of all graduate mathematics programs. NSF HERD survey data indicate that if only federal sources are included, as with WVU, more than 50\% of R1 mathematics departments would be under similar reviews and in jeopardy. SERD funding levels contribute to the rankings as described in~\cite{ams2023sm} and in particular to the Aggregated Research Index Score (ARIS) used as one factor for R1 status since it is compiled from external funding, postdoctoral researchers and PhDs conferred. This index and the rank of an IHE increase as external funding grows, but the contributions from mathematics departments and graduate programs in mathematics, however, extend well beyond research funding, and many other factors derived from a department also contribute to ranking and reputation. Figure \ref{fig:matharisTT}, for example, shows that the number of TT faculty in mathematics by itself correlates positively with increases in the ARIS aggregate measure of research productivity in the R1 grouping.

Similar analysis~\cite{ams2023sm} shows that other factors such as lower student-faculty ratio in undergraduate mathematics courses, higher numbers of graduate students in mathematics, and higher numbers of mathematics PhDs produced also positively correlate with ARIS ranking, supporting the notion that an active and well-resourced research program in mathematics supports R1 status in ways that are independent of research funding levels. These correlations also give some insight into why 97\% (142 of 146) of R1 IHEs maintained graduate programs conferring doctorates in mathematics up until $2023$. Students and other researchers benefit from access to advanced mathematics, and this access enhances research programs in multiple STEM departments, economics and business as well as health sciences. 

In order to meet their instructional needs in mathematics, mathematics departments within R1 IHEs appoint a combination of tenured or tenure-track faculty (TT) with non tenure-track teaching faculty (NTT). At most of the IHEs in the United States, departments engage in the instruction of service courses such as calculus, differential equations, or linear algebra to sometimes thousands or tens of thousands of students per year and this demands a strategic balance of faculty with different teaching loads and research expectations. The populations of TT and NTT mathematics faculty along with the ratio of undergraduate students in all math courses to total math faculty (UG/Math FAC) are shown in Table \ref{tab:mathr1}.  Other metrics such as the ratio of TT to NTT faculty, average federal research funding (RF, in \$1000s) in mathematics, the mean number of postdoctoral researchers in math (Mean PD), the five-year mean number of mathematics PhDs conferred (Mean PhDs), UG mathematics majors (UG math \%) as a percent of the total population of UG students taking math courses, and graduate PhD and MS mathematics students (GR math \%) as a percentage of all students in graduate mathematics courses are shown as well.
\begin{table}[h!]
\centering
\footnotesize
\begin{tabular}{p{.22\columnwidth}cccc}
& \thead{\small Private} & \thead{\small Public} & \thead{\small HSI}   & \thead{\small LG} \\
 UG/Math FAC& 61.5    & 131.8     & 131.2 & 124.5        \\
 Math RF& 2927   & 1948.7     & 727.1 & 2748.6        \\
 Mean TT               & 28.7    & 37.9      & 44.6  & 43.4         \\
 Mean NTT              & 7.4     & 15.3      & 15.6  & 16.7         \\
 TT/NTT            & 5.5     & 3.9       & 4.3   & 5.0          \\
 Mean PD         & 13.1    & 7.5       & 4.4   & 10.6          \\
 Mean PhDs       & 10.0    & 10.6      & 9.1   & 12.3          \\
 UG math \%& 31.5\%    & 30.9\%    & 25.7\% & 32.2\%        \\
 GR math \%& 4.6\%    & 5.2\%      & 10.3\% & 4.8\%      
\end{tabular}
\caption{Characteristics of mathematics departments at R1s}
\label{tab:mathr1}
\end{table}
For land-grant IHEs, the ratio of undergraduates in math courses to total department faculty in the first row of Table \ref{tab:mathr1} correlates negatively with the ARIS research ranking of the university for values greater than 150:1~\cite{ams2023sm}. As the ratio increases, aggregate research indicators begin to drop, suggesting that an institution with ratios nearing 200:1 may have structural issues with the ability of research faculty to participate in research projects and to procure external funding. It is unclear how effectively an R1 can then lower the number of both TT and NTT faculty in mathematics as WVU plans to do and still serve its educational community properly or even meet the instructional needs of its student body in mathematics. With ratios nearing or above 300:1, all STEM disciplines begin to struggle with insufficient access to needed mathematics coursework for their students, and the profile of the university as a whole changes. R1 status depends, at least in part, on the trends shown here for mathematics resources. The way in which R1 mathematics departments meet instructional needs interacts synergistically with measures of research productivity and these include both the number of TT faculty in a department as well as the numbers of graduate students, PhDs conferred and grant funding. The actions proposed at the IHEs listed in the introduction, however, and in particular at WVU, indicate policy choices that assume that R1 IHEs can operate without core graduate and research capacities in mathematics as well as other areas. The motivations for these changes reflect more than funding shortfalls, and ultimately indicate a belief that portions of those departmental missions can be sacrificed without damaging others. Increases in R1 ranking factors can be obtained from other sources of external funding such as health sciences, and, as noted in the previous section, state funding can be realigned to support the RAM portions of a university even a state system. This appears to be the case for WV. The actions taken in response to recent funding concerns within land-grants such as WVU conflict with their core missions primarily because they limit student access to the best undergraduate course opportunities, but they also eliminate access to graduate study and potentially limit research capacity overall. Ignoring the role of mathematics as discussed above becomes especially problematic in the context of land-grant R1s where the research and teaching missions exist primarily to provide the best opportunities to the citizens of a state. Relying on R1 status alone as an indicator of the excellence in academic programs, or as a basis for maintaining any program at an IHE, then becomes problematic since that status may only reflect research expenditures that are completely independent of mathematics, physics, engineering or any scholarly work within those programs. For land-grant, HSI, HBCUs, and other public IHEs with much broader missions, it is essential that resource allocation take all missions into account, and that measures of excellence incorporate multiple factors including grant supported research programs as well as a variety of other indicators. Departments with research programs must continue to seek funding for their research, but more generally departments should work to increase the appeal of mathematics as a major, incorporate more computational and data science components into mathematics programs, and work to retain majors at critical departure points.  Penalizing departments with program discontinuation and terminating productive and talented faculty, however, punishes students by limiting their options unnecessarily. 

\section*{Conclusions}
Internal or external forces may change resource allocations and departments must become more adept at working to demonstrate impact in order to have a voice in those allocations. A vice president for research can reallocate the funding for a graduate program to another department or a provost can redirect faculty resources away from STEM programs and into other areas in order to meet what they feel is a better strategic goal for a university, and mathematics departments must be ready to engage with other departments and administrators to demonstrate the implications of those choices. We must also clearly articulate outside of the university in legislative circles and with the public the value of the contributions of mathematics to education and research in general. We must be ready to provide evidence that supports the effectiveness of our teaching as well as outcomes that demonstrate the impact of our research and the value of our degrees~\cite{ams2023sm}. Without this dialogue, we risk being characterized as expendable. At the core of large public institutions, our mission is to serve the people of our communities by providing opportunities to study and learn in areas that will benefit society in ways that sometimes are not reflected in the number of people that graduate from a specific degree program or the total amount of grant funding supporting the research of a department. Narratives that focus on a single aspect of a department or ignore significant portions of what a department does will not effectively support the missions of an IHE.

The administrators at WVU and the state leaders of West Virginia overlooked the contributions of mathematicians to the educational opportunities of the people in that state~\cite{ams2023sm} and to its flagship university. It may seem at first glance that this is a local concern dependent on financial or enrollment issues, but these arguments appear in other states and localities, spread in part by consultants with no commitment to the communities those universities serve. External funding from grants or federal appropriations are important for the support of research and even graduate programs, but the pursuit of that funding for its own sake may reinforce alarming shifts away from academics by creating the notion that all efforts must be aligned with that pursuit. The metrics in that case become the goal~\cite{ams2023sm} of those operations with no consideration for the quality of core academics. In the example of WVU, leaders have aligned their choices with health science operations, medical school support, and with degrees they perceive to be the most employable options for students. These supposed efforts to ``modernize" programs or provide students with the ``best" career options become arguments that ignore the essence of what an IHE does. How academia responds may determine which universities continue to be institutions of true learning and which become minimally effective professional training programs. It is possible to build a ``Very High Activity" research institution around a health sciences paradigm or industrial applications, but in the case of at least public land-grant universities, it ignores the fundamental intent of their creation. Gordon Gee, president of WVU, responded~\cite{ams2023sm} to concerns about discontinuing the PhD in Mathematics by saying ``Someone else is going to have a great PhD program in mathematics. And you know what? God bless them." True flagship universities with effective leadership see opportunity in providing the most \textit{liberal and practical education} possible and move towards it, not away. The students of states like West Virginia will have to look elsewhere to find that education, and for those students the lack of access is a failure that is not reflected in any financial analysis, research status, or ranking.  The ``academic transformation" of WVU is a model for changes to the underlying strategic mission of IHEs that does not reflect the original intent of publicly funded higher education as spending in essential areas such as mathematics, world languages or libraries are further reduced to compensate for interest payments on growing debt or private aircraft travel for administrators~\cite{wvuair2023}. 

It remains to be seen what the long term impact of the program cuts and faculty reductions will be at WVU. The overall structure of large, public R1s continues to attract students as shown in Figure \ref{fig:enroll} for academic as well as non-academic reasons. IHEs such as WVU can currently maintain R1 status and emphasize other programs, at least in principle, but the students and faculty will ultimately decide if the quality of what remains meets their needs or if they will search elsewhere for academic opportunity. The guidelines for Carnegie classification are being revised to include ``a wider set of dimensions that define an institution, and not distort the process by overweighting one single aspect of an institution’s purpose and activities"~\cite{newR12023} so it may be that a university that focuses solely on the value of grants obtained may be classified as that, and one with a broad collection of graduate offerings will be recognized that way. Removing graduate mathematics redefines WVU, and it is arguable that the loss of mathematics as well as world languages and other fields diminishes the university in important ways. The outcomes at other IHEs such as Emporia State University in Kansas where student enrollments dropped another $12.5$\%~\cite{ams2023sm} after program changes suggest that universities considering similar actions should do so carefully. Mathematically speaking, they may accelerate the declines to which they are responding instead of turning them around.

Just as the value of an education is more than the credential acquired, the value of mathematics extends from more than just mathematics course requirements in a degree program or even a degree in mathematics. Mathematics gives students an opportunity to understand ideas ranging from differential equations to graph theory and to use these ideas to solve problems in finance, renewable energy technology or molecular structure prediction. Mathematics is necessary for IHEs to make those educational opportunities available, and other disciplines need advanced mathematics to facilitate their research. It exposes students to a wealth of knowledge that begins with mathematics but spreads out across all areas of scientific inquiry. Planning and resource allocation at an IHE must account for both the core educational impact of mathematics as well as the broader, sometimes harder to measure, impact of mathematics research on other disciplines and student opportunities in order for an institution to excel. 

\subsection*{Acknowledgements}
Many thanks to Cristian Potter for his help in gathering the data related to mathematics departments at R1s, and to Karen Saxe for helpful suggestions and comments during the writing of this article. Additional thanks are extended to the many investigative reporters, journalists, faculty and students who have published the many works documenting different aspects of the situation at WVU collected more completely here~\cite{ams2023sm}. 

\section*{Additional Discusson}
\subsection*{Mathematics Impacts}
As noted in ~\cite{fullerNotices24}, most departments of mathematics have substantial teaching missions that are derived from the foundational role that mathematics plays in almost every degree program. Ultimately, some aspect of mathematics sits at the core of the essential workings of most careers even if computing derivatives of trigonometric functions doesn't occur regularly in the workplace. Reasoning and critical thinking are necessary skills, and mathematics serves as one of the strongest areas within which a student can develop them. Due in part to the prominence of coursework in degree programs and in part to the complex nature of the subject, mathematics departments in higher education also have a history of being a barrier to student success~\cite{olson2012engage}. This can lead to conflict that focuses on student success in our courses, but the last ten years or so have seen significant work leading to knowledge about best practices~\cite{rasmussen2019brief,bressoud2013calculus} for mathematics education as well as implementations that lead to change~\cite{smith2017seminal,calcscience}. These outcomes build on substantial prior work~\cite{hurley1999effects} giving an indication that the profession has inertia in this area. Mathematics research also contributes in a multitude of areas~\cite{andjelkovic2020topology}, and these impacts lead to ever widening options for mathematics-related research, careers~\cite{cass2011examining,leyva2022making,usnewsmath} and application paths.\cite{menezes2018handbook,baz2004financial,axioms12070699}. Research in mathematic ranges from areas that exist as the development of new theories to applications that exist primarily to solve problems in other disciplines. Taken together, these impact make the discipline one of the most critical aspects of higher as well as K-12 education and the infrastructure of almost any educational enterprise. 


\subsection*{The Events at West Virginia University (WVU)}
The recent review of the Department of Mathematics at West Virginia University and of other departments at other universities such as the University of North Carolina-Greensboro have become focal points in discussions surrounding budget driven restructuring of large, public, land-grant IHEs. I served as chair of the Department of Mathematics at WVU from 2008-2018, and led a number of initiatives in research, graduate education, and external funding. The narrative being reported in media outlets and by the university criticizes the now School of Mathematics and Data Science for shortfalls, but in many cases the underlying strategic actions of the unit supported all the outcomes sought by the administration and the consultancies reviewing the university. Performance reviews occur at most if not all IHEs, and when designed with the intent to build and create excellence, program reviews often yield helpful insights that can be used to improve a department and even a university. When designed without strategic vision, the disconnect between performance goals that support the core mission of higher education or even the specific core missions of an individual IHE and the goals used to evaluate departments can be quite substantial. This leads to conflicting approaches to those missions and to undesirable outcomes.



Over the period of time since its first PhD conferral in $1993$, the Department of Mathematics at WVU worked steadily towards building robust research capacity and had become one of the strongest graph theory departments in the country if not the world. Many of its researchers were internationally involved and research productivity had grown so that in $2014$ the graduate program was ranked nationally for the first time in its history. Publication in peer-reviewed journals had grown by more than 30\% since 2000, and total research expenditures grew from \$460,000 in 2013 at the beginning of my second term to over \$1m in 2016 and more than \$2.5m in 2018 through strategic investment in faculty hiring, pre-award support, and seed grant programs. At the same time, graduate enrollment grew form 41 students to 96 in 2017 at the maximum of the program.  

\begin{figure}[h!]
\centering
\includegraphics[width=.9\columnwidth]{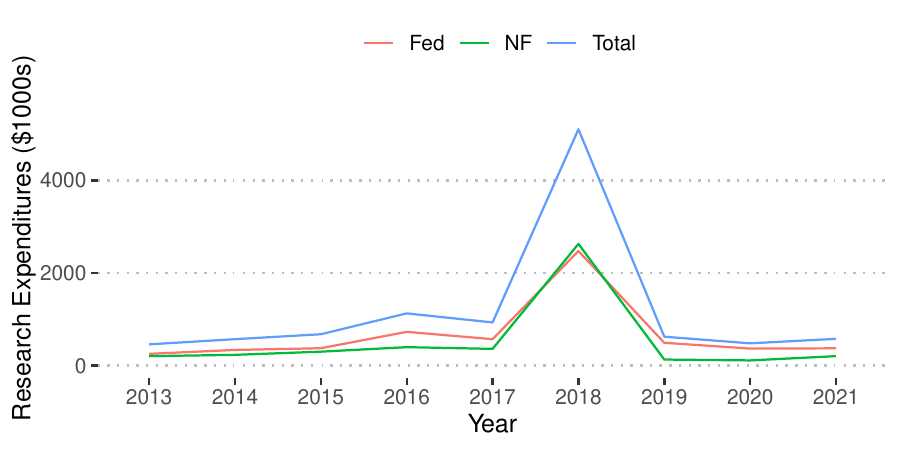}
\caption{Total Grant Funded Research Expenditures at WVU from 2011 to 2021}\label{fig:wvuserd}
\end{figure}

\noindent Some of that capacity was lost due to faculty departures. Strategic changes at the university level led to difficulty maintaining programs, and the onset of Covid-19 contributed to enrollment struggles. 

It is unclear if that capacity will return in the current university environment. Many articles related to the program reductions at WVU have appeared in various news outlets including the New York Times~\cite{nytwvu} and the Boston Review~\cite{bosrev2023} . The Chronicle of Higher Education has published multiple pieces and opinions including~\cite{chron2015}. The specifics of the program reductions can be found in many of these as well as~\cite{wvurec} and~\cite{quinn2023}. Arkansas and Missouri have undergone similar reviews and in some cases reduced programs~\cite{hend2022}.  Proposals to reduce the degree requirements in mathematics are discussed in~\cite{caks2022} and further program reduction recommendations are found in~\cite{mckuf}. The impact of these reductions at WVU are discussed here~\cite{fightwvu2023} and the responses from WVU President E. Gordon Gee have been published here~\cite{Krebs2023} and here~\cite{timeswv2023}. Other universities are experiencing similar reviews, and many are overseen by rpkGroup, the consulting firm retained by WVU. UNC-Greensboro has developed recommendations to discontinue multiple programs including Physics~\cite{uncg2023}. In the middle of these program review debates, the state of WV has passed laws revising the tenure structure at WVU~\cite{tenurewvu2022} in order to allow faculty to be terminated. These changes are being replicated in different ways in other states including Wisconsin~\cite{tenure2017}, Texas~\cite{tenuretx2023} and now Florida.

\subsection*{Carnegie Classifications and R1 Status}
The Carnegie Classification System was developed in 1973 by the Carnegie Commission on Higher Education but has been administered by different organizations~\cite{AceNetCC} since then.  A total of $146$ IHEs of the $3939$ in the US in $2021$ were classified as R1 by the Indiana University Center for Postsecondary Research, or about $3.7\%$. The evaluation of IHEs begins with the conferral of at least $25$ doctorates, and then ranks those IHEs using the number of faculty in ladder-rank appointments, the number of postdoctoral researchers (PD), the amount of external funding, and the number of PhDs conferred. Differences in the characteristics of the institutions are shown in Table \ref{tab:R1char} where they are grouped into private, public land-grant, non-land grant HSIs, and other public IHEs. The table summarizes the mean values of several key aspects of R1 IHEs including total undergraduate enrollment (UG), total graduate enrollment (GR), Bachelors (BS) and PhDs (PhD) awarded, and the number of full-time faculty in ladder rank (assistant, associate, and full professors). The mean number of students enrolled at the undergraduate level provides the clearest discrepancy between the public and private R1s. External funding listed here as Science and Engineering Research and Development Expenditures  (SERD, in millions of dollars) are considerably higher for private and land-grants, while the average number of faculty and the number of PhDs awarded are almost the same. This metric along with the number of postdoctoral non-faculty researchers (PD) is the largest factor impacting R1 ranking as of the $2021$ analysis.

\begin{table}[h!]
\centering
\footnotesize
\begin{tabular}{p{.2\columnwidth}cccc}
& \thead{\small Private} & \thead{\small Public} & \thead{\small HSI} & \thead{\small LG} \\
 Number       & 39     & 107      & 11       & 51      \\
 SERD   & 579.7  & 388.8   & 235.6   & 483.7   \\
 UG& 8966 & 25241 & 28382 & 26331 \\
 GR& 9315 & 7787  & 6838  & 7792  \\
 Mean PhD     & 331  & 345   & 287   & 415   \\
 Mean BS      & 2629 & 6228  & 6882  & 6700  \\
 Mean FAC     & 1471 & 1474  & 962   & 1624  \\
 Mean PD      & 887  & 443  & 299   & 580   \\
\end{tabular}
\caption{Characteristics of R1 Universities by Group with SERD Research Funding in Millions, all Undergraduate (UG) and Graduate (GR) enrollments, mean Ladder Rank Faculty (FAC), degrees granted, and Postdocs (PD).}
\label{tab:R1char}
\end{table}

\begin{figure}[h!]
\centering
\includegraphics[width=.9\columnwidth]{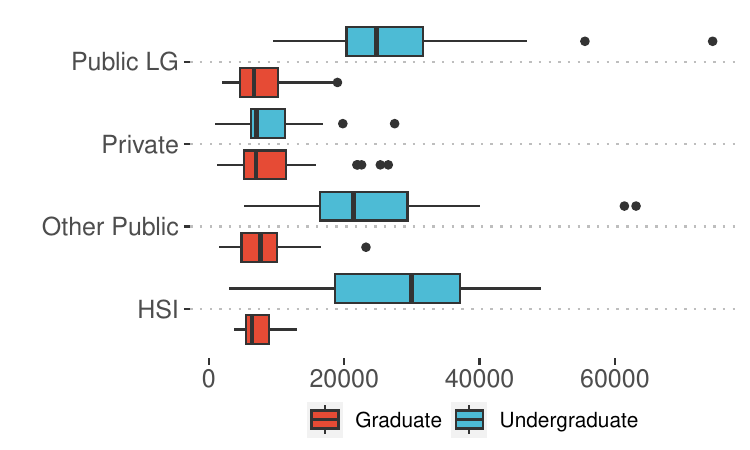}
\caption{Undergraduate and Graduate Enrollments at R1 IHEs by Institution Type}\label{fig:ugenr}
\end{figure}

\subsection*{The ARIS and PCIS Weighted Rank Scores}
\begin{figure}[h!]
\centering
\includegraphics[width=.9\columnwidth]{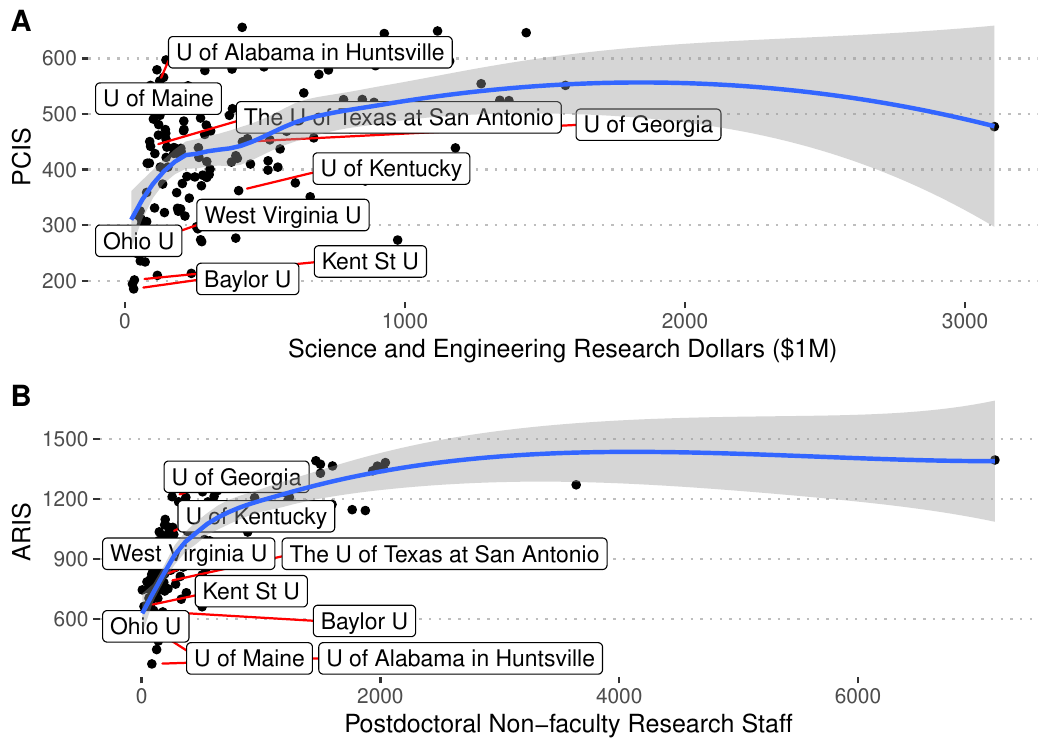}
\caption{a)PCIS by Science and Engineering Research and Development Expenditures (millions) and b)ARIS by Postdoctoral Non-Faculty Researchers}\label{fig:serd}
\end{figure}

Most of these factors, however, do not impact R1 ranking in a substantive way. R1 status in the $2021$ rankings was determined by a scaled Euclidean distance computed from two standardized, weighted rank scores: the Per Capita Index Score (PCIS) and the Aggregate Research Index Score (ARIS). PCIS uses external funding and postdoctoral research numbers per ladder rank faculty to create a weighted rank score, and ARIS creates a weighted composite rank of raw external funding, postdoctoral researcher numbers, and doctorates granted in multiple areas. Figure \ref{fig:serd} shows the relationship between the PCIS factor and research funding along with a plot of ARIS against the number of postdoctoral researchers. Since the horizontal values are included in both vertical axis factors, both of these inputs impact the R1 ranking process. As a result, many universities invest heavily in procuring external funding and non-faculty researchers to ensure high R1 rankings.

There are other correlative relationships, some of which are shown below, but in many cases changes in departmental numbers or college level totals do not impact R1 rankings significantly. In the case of a university like WVU, for example, reducing the total college of arts and sciences research funding to zero and removing all conferred STEM doctorates would only lower the overall rank score from $0.35$ to $-0.025$. The lower bound for R1s is set at -0.15 in the 2021 rankings. Indeed, the research profile of IHEs varies widely, with many of them focusing on health science related expenditures or other industry support. The focus on these factors, however, can translate into resource reallocation at the cost of core academics. In some sense, the focus on either external funding or non-faculty postdoctoral researchers alone is an illustration of Goodhart's Law~\cite{chrystal2003goodhart} where the measure becomes the goal at the expense of other important factors.

Looking at the ARIS index plotted against the average total research funding levels over the last three years for R1 mathematics departments shows that increasing departmental external funding levels positively correlates with increases in these rankings. 
\begin{figure}[h!]
\centering
\includegraphics[width=.9\columnwidth]{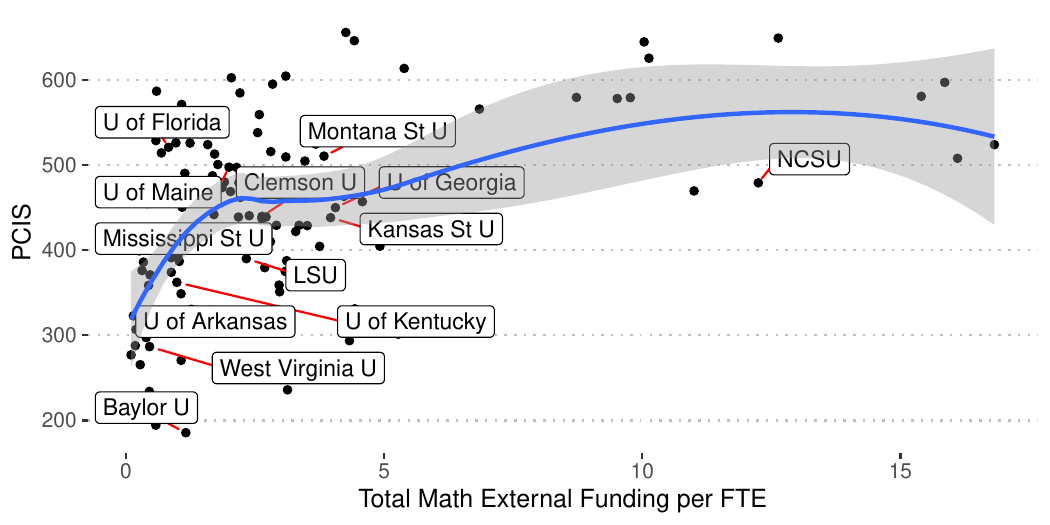}
\caption{Scatterplot of the PCIS Weighted Ranking by Mathematics Research Funding per University FTE}\label{fig:mathpcis}
\end{figure}

\begin{figure}[h!]
\centering
\includegraphics[width=.9\columnwidth]{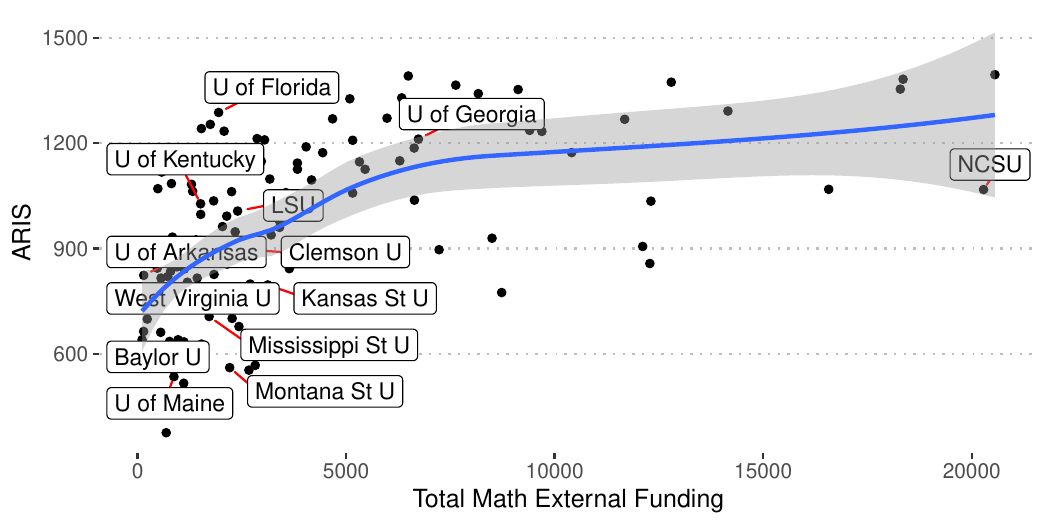}
\caption{Scatterplot of the ARIS Weighted Ranking by Mathematics Research Funding}\label{fig:matharisserd}
\end{figure}

\subsection*{Goodhart's Law and Chasing Metrics}
When a specific metric is identified as measuring a desired result, over time it can become the case that evaluations begin to only reflect the value of the metric and not the underlying state that the metric was initially chosen to represent~\cite{chrystal2003goodhart}. Research spending at R1s consitutes a good indicator of research activity, but the number of peer-reviewed publications remains a core measure of the output of the research enterprise including citations, presentations and publications in conferences, and other more difficult to measure interactions with colleagues and funding agencies. Establishing the value of research funding as a sole metric of success may lead to situations where investments of resources are only made that lead to increases in those numbers, and that investment may or may not lead to true research outcomes.

Similarly, viewing student success as a percentage of grades of A, B or C alone also becomes a self-fulfilling effort that may not reflect the substance of the learning experience for students. Learning outcomes are complex, nuanced and difficult to measure, and student experience is reflected by more than just student evaluations of instruction. Changing instructional culture solely to increase student grades may not lead to improvements in learning, and can lead to weaker programs overall.

\subsection*{Additional Faculty Contributions}
Within this context, more complex considerations have begun to evolve concerning the contributions that a discipline like mathematics makes to the mission of a university. The recent events noted above raise concerns about changing perceptions of the role of a mathematics department. When leaders at flagship state IHEs such as President Ben Sasse~\cite{sassenyt2023} at the University of Florida question ask 
\begin{quote}
\it What will today’s generic term ‘professor’ mean when you disaggregate syllabus designer, sage-on-the-stage lecturer, seminar leader, instructional technologist, grader, assessor, etc.?
\end{quote}
implies many things but it takes pieces of faculty activities and reduces them to simplistic notions that hide the true difficulty and craft in the processes.
We can ask what value a mathematics department and its various collective skillsets brings to a university and compare some to the propensity of an IHE to participate in very high levels of research activity. We can also gain insight into whether or not they sum to more than the list of parts quoted above.
 
\begin{table}[h!]
\centering
\footnotesize
\begin{tabular}{p{.22\columnwidth}cccc}
& \thead{\small Private} & \thead{\small Public\\ LG} & \thead{\small HSI}   & \thead{\small Other\\ Public} \\
 UG/Math FAC& 61.5    & 128.1     & 146.2 & 133.9        \\
 Math RF& 2927   & 2682     & 637 & 1468        \\
 Mean TT               & 28.7    & 43.3      & 45.8  & 32.3         \\
 Mean NTT              & 7.4     & 17.3      & 17.0  & 13.5         \\
 TT/NTT            & 5.5     & 4.7       & 1.5   & 3.3          \\
 Mean PD         & 13.1    & 9.3       & 0.5   & 6.5          \\
 Mean PhDs       & 10.0    & 11.9      & 6.5   & 9.7          \\
 UG major \%& 32\%    & 32\%    & 25\% & 31\%        \\
 GR major \%& 5\%    & 5\%      & 15\% & 5\%      
\end{tabular}
\caption{Mean Characteristics of Math Units in R1s including the Average Number of Undergraduates (UG) in Math Courses per Math Faculty (FAC, all ranks and types)}
\label{tab:mathr1}
\end{table}

For land-grants such as WVU, the ratio of undergraduates in math courses to faculty (the first line in Table \ref{tab:mathr1}) correlates strongly with the research profile of the university. As the ratio increases, aggregate research indicators begin to drop, indicating unbalanced resource allocation. An institution with a ratio as high as that of WVU, near 200:1, may be over-stressing faculty and their ability to participate in research projects in the department and otherwise. It is unclear how effectively an R1 can then lower the number of both TT and NTT faculty in those roles as WVU plans and still serve its educational community properly or even meet the instructional needs of its student body in mathematics. With ratios nearing or above 300:1, all STEM disciplines begin to struggle with insufficient access to needed mathematics coursework for their students, and the profile of the university as a whole changes.

\subsection*{Student Share of Costs}
\begin{figure}[h!]
\centering
\includegraphics[width=.9\columnwidth]{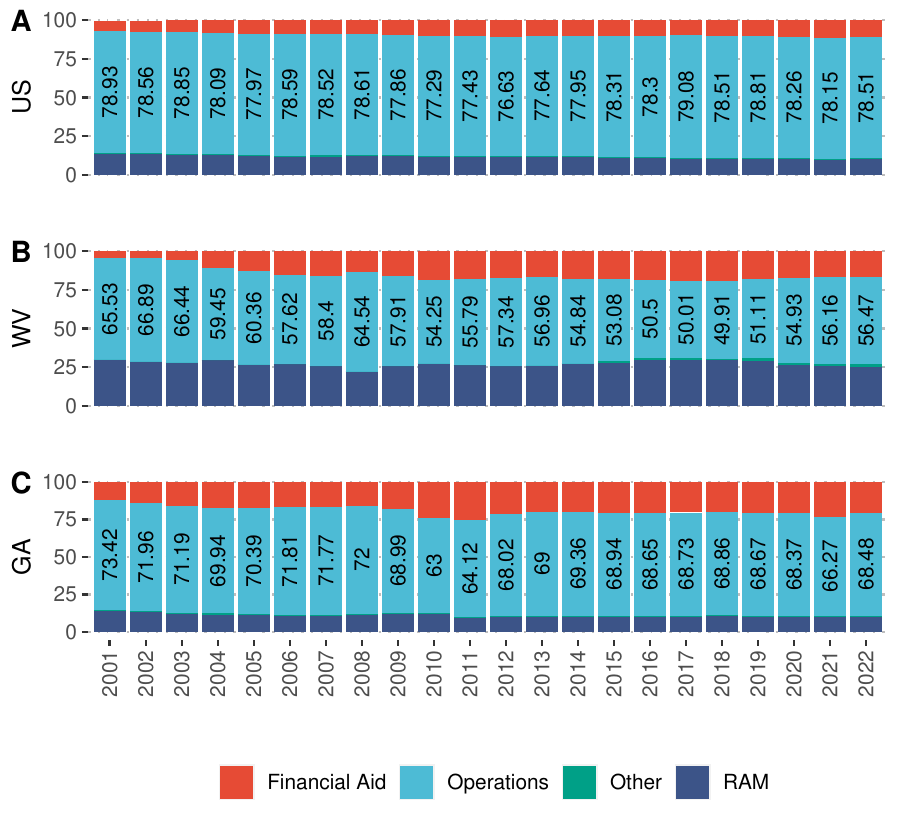}
\caption{Distribution of State Appropriated Funds in the a) US Overall, b) West Virginia and c) Georgia 2001 - 2022}\label{fig:distfunds}
\end{figure}

Public funding for IHEs has been decreasing~\cite{pbs2019} for many decades, and it has been argued that many of the resource conflicts that even large, public IHEs are beginning to experience are the result of many years of policy changes.  As shown in Figure \ref{fig:distfunds}, state and university leaders have dramatically shifted support in ways that are beginning to reduce operational capacity at even large institutions. 

\begin{figure}[h!]
\centering
\begin{subfigure}{\columnwidth}
\includegraphics[width=.9\columnwidth]{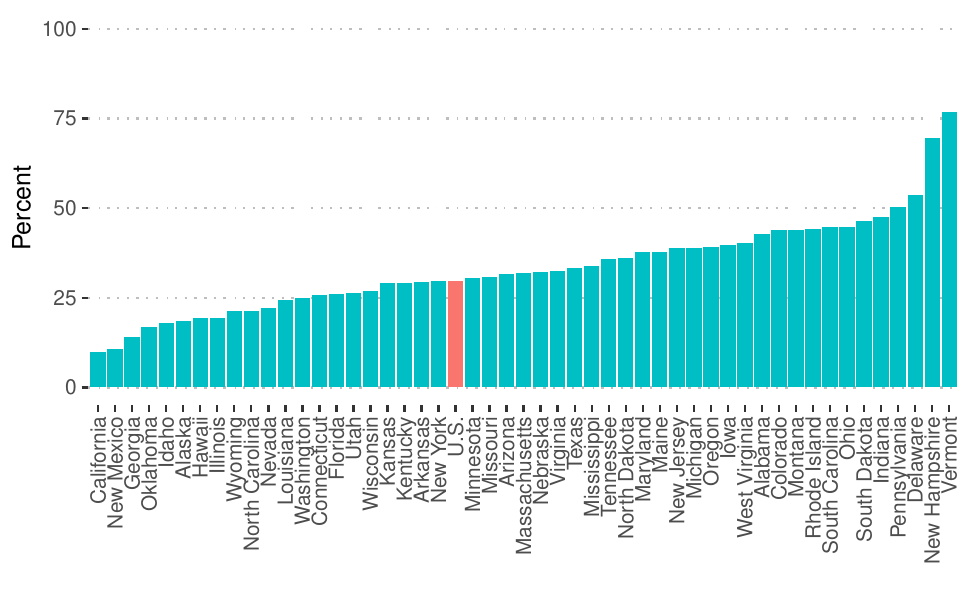}
\caption{Student Share by State 2002}\label{fig:sscost02}
\end{subfigure}
\begin{subfigure}{\columnwidth}
\includegraphics[width=.9\columnwidth]{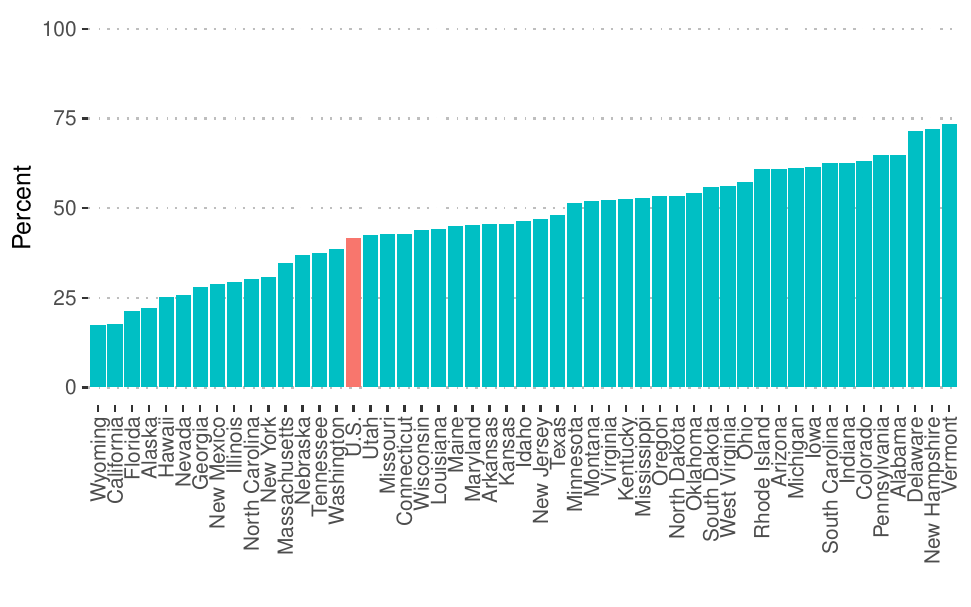}
\caption{Student Share by State 2022}\label{fig:sscost22}
\end{subfigure}
\caption{Student Share of State Education Revenue - Net Tuition as a Fraction of Total State Education Revenue}
\end{figure}

This has led to dramatic increases in student tuition, and as these student costs increase, ironically student access to educational options are dramatically reduced. The restructuring in WV~\cite{unkept2018} is perhaps one of the most dramatic examples of this effect, where legislative appropriations coupled with resource shifts have led to program changes across the entire university. The enrollment shifts that WVU and other IHEs have seen certainly reflect some part of the demographic shift that has been characterized as the 'enrollment cliff', but some part of it relates to the changing access to programs and educational opportunities there. 

\subsection*{Data Sources}
Data used in these analyses can be obtained from the NSF HERD survey~\cite{herd2023} and IPEDS data\cite{ipeds} repository. Additional information related to higher education funding can be found in the SHEFF report~\cite{sheff2023}. Data related to mathematics department staffing, student levels, and more detailed information was obtained from the Conference Board on the Mathematical Sciences and American Mathematical Society annual survey~\cite{cbmsdata}.

\bibliographystyle{plain}
\bibliography{References}

\begin{bibdiv}
\begin{biblist}

\bib{hend2022}{misc}{
      author={Adame, Jaime},
      author={Earley, Neal},
       title={{ASU System approves Henderson State cuts; faculty jobs, degree
  programs slashed}},
   publisher={Arkansas Democrat-Gazette},
        date={2022},
  url={https://www.arkansasonline.com/news/2022/may/06/asu-system-approves-henderson-state-cuts/},
  note={\url{https://www.arkansasonline.com/news/2022/may/06/asu-system-approves-henderson-state-cuts/},
  Accessed October 4, 2023},
}

\bib{morrill1862}{misc}{
      author={Archives, National},
       title={{Morrill Act} of 1862},
        date={1862},
         url={\url{https://www.archives.gov/milestone-documents/morrill-act}},
        note={\url{https://www.archives.gov/milestone-documents/morrill-act},
  Accessed September 5, 2023},
}

\bib{andjelkovic2020topology}{article}{
      author={Andjelkovi{\'c}, Miroslav},
      author={Tadi{\'c}, Bosiljka},
      author={Melnik, Roderick},
       title={The topology of higher-order complexes associated with brain hubs
  in human connectomes},
        date={2020},
     journal={Scientific reports},
      volume={10},
      number={1},
       pages={17320},
}

\bib{baz2004financial}{book}{
      author={Baz, Jamil},
      author={Chacko, George},
       title={Financial derivatives: pricing, applications, and mathematics},
   publisher={Cambridge University Press},
        date={2004},
}

\bib{bressoud2013calculus}{article}{
      author={Bressoud, David~M},
      author={Carlson, Marilyn~P},
      author={Mesa, Vilma},
      author={Rasmussen, Chris},
       title={The calculus student: insights from the {Mathematical Association
  of America} national study},
        date={2013},
     journal={International Journal of Mathematical Education in Science and
  Technology},
      volume={44},
      number={5},
       pages={685\ndash 698},
}

\bib{bressoud2012second}{book}{
      author={Bressoud, David~M},
       title={Second year calculus: from celestial mechanics to special
  relativity},
   publisher={Springer Science \& Business Media},
        date={2012},
}

\bib{forbes2023}{misc}{
      author={Bauer-Wolf, Jeremy},
       title={{University of Kansas} looks to cut 42 academic programs},
   publisher={Higher Ed Dive},
        date={2023},
  url={https://www.highereddive.com/news/16-states-underfunded-land-grant-hbcus-by-over-12b-biden-admin-says/694020/},
  note={\url{https://www.highereddive.com/news/16-states-underfunded-land-grant-hbcus-by-over-12b-biden-admin-says/694020/},
  Accessed September 19, 2023},
}

\bib{cbmsdata}{misc}{
      author={CBMS},
       title={Conference board of mathematical sciences survey},
        date={2024},
        note={Retrieved April 2, 2024 from
  \url{https://www.ams.org/profession/data/cbms-survey/cbms-survey}},
}

\bib{cass2011examining}{inproceedings}{
      author={Cass, Cheryl~AP},
      author={Hazari, Zahra},
      author={Cribbs, Jennifer},
      author={Sadler, Philip~M},
      author={Sonnert, Gerhard},
       title={Examining the impact of mathematics identity on the choice of
  engineering careers for male and female students},
organization={IEEE},
        date={2011},
   booktitle={2011 frontiers in education conference (fie)},
       pages={F2H\ndash 1},
}

\bib{chrystal2003goodhart}{article}{
      author={Chrystal, K~Alec},
      author={Mizen, Paul~D},
      author={Mizen, PD},
       title={{Goodhart’s Law}: its origins, meaning and implications for
  monetary policy},
        date={2003},
     journal={Central banking, monetary theory and practice: Essays in honour
  of Charles Goodhart},
      volume={1},
       pages={221\ndash 243},
}

\bib{bosrev2023}{misc}{
      author={Casey, Rose},
      author={Wilkerson, Jessica},
      author={Winant, Johanna},
       title={An open letter from faculty at {West Virginia University}},
   publisher={Boston Review},
        date={2023},
  url={\url{https://www.bostonreview.net/articles/an-open-letter-from-faculty-at-west-virginia-university/}},
  note={\url{https://www.bostonreview.net/articles/an-open-letter-from-faculty-at-west-virginia-university/},
  Accessed September 8, 2023},
}

\bib{tenuretx2023}{misc}{
      author={Ebbeler, Jennifer},
       title={{What Ending Tenure Would Mean to Texas}},
   publisher={Inside Higher Ed},
        date={2023},
  url={https://www.insidehighered.com/opinion/views/2023/04/20/what-ending-tenure-would-mean-texas},
  note={\url{https://www.insidehighered.com/opinion/views/2023/04/20/what-ending-tenure-would-mean-texas},
  Accessed October 4, 2023},
}

\bib{timeswv2023}{misc}{
      author={Fernandez, Esteban},
       title={{WVU Math Department appeal chaos drives wedge between faculty
  and administrators}},
   publisher={Times West Virginian},
        date={2023},
         url={https://timeswv-cnhi.newsmemory.com/?publink=190472773_134ad87},
  note={\url{https://timeswv-cnhi.newsmemory.com/?publink=190472773_134ad87},
  Accessed October 12, 2023},
}

\bib{ipeds}{misc}{
      author={for Education~Statistics, National~Center},
       title={Integrated postsecondary education data system {(IPEDS)}},
        date={2023},
         url={\url{https://nces.ed.gov/ipeds/}},
        note={Accessed September 4, 2023},
}

\bib{tenure2017}{misc}{
      author={Flaherty, Colleen},
       title={{Killing Tenure}},
   publisher={Inside Higher Ed},
        date={2017},
  url={https://www.insidehighered.com/news/2017/01/13/legislation-two-states-seeks-eliminate-tenure-public-higher-education},
  note={\url{https://www.insidehighered.com/news/2017/01/13/legislation-two-states-seeks-eliminate-tenure-public-higher-education},
  Accessed October 4, 2023},
}

\bib{tenurewvu2022}{misc}{
      author={Flaherty, Colleen},
       title={{Shoring Up Tenure, or Weakening It?}},
   publisher={Inside Higher Ed},
        date={2017},
  url={https://www.insidehighered.com/news/2022/10/27/west-virginia-u-proposal-outlines-process-firing-faculty},
  note={\url{https://www.insidehighered.com/news/2022/10/27/west-virginia-u-proposal-outlines-process-firing-faculty},
  Accessed October 4, 2023},
}

\bib{herd2023}{misc}{
      author={for Science, National~Center},
      author={Statistics, Engineering},
       title={Higher education research and development {(HERD)} survey},
   publisher={National Science Foundation},
        date={2023},
  url={https://ncses.nsf.gov/surveys/higher-education-research-development/2021},
  note={\url{https://ncses.nsf.gov/surveys/higher-education-research-development/2021},
  Accessed September 24, 2023},
}

\bib{fullerNotices24}{article}{
      author={Fuller, Edgar},
       title={Are university budget cuts becoming a threat to mathematics?},
        date={2024May},
     journal={Notices of the American Mathematical Society},
      volume={71},
      number={5},
       pages={TBD},
        note={to appear, April 2024},
}

\bib{ams2023sm}{misc}{
      author={Fuller, Edgar~J.},
       title={{ Are Math Departments At Risk During University Budget Cuts?
  Extended Version}},
        date={2024},
}

\bib{garcia2019defining}{article}{
      author={Garcia, Gina~A},
       title={Defining “servingness” at {Hispanic}-serving institutions
  ({HSIs}): Practical implications for {HSI} leaders},
        date={2019},
     journal={American Council on Education},
}

\bib{hurley1999effects}{article}{
      author={Hurley, James~F},
      author={Koehn, Uwe},
      author={Ganter, Susan~L},
       title={Effects of calculus reform: Local and national},
        date={1999},
     journal={The American Mathematical Monthly},
      volume={106},
      number={9},
       pages={800\ndash 811},
}

\bib{jackson1997whatever}{article}{
      author={Jackson, Allyn},
       title={{Whatever happened to Rochester? Two years later, mathematics is
  getting accolades}},
        date={1997},
     journal={Notices of the American Mathematical Society},
      volume={44},
      number={11},
       pages={1463\ndash 1466},
}

\bib{rochester1996}{misc}{
      author={Jaffe, Arthur},
      author={Lipman, Joseph},
      author={},
      author={Lowengrub, Morton},
       title={{U. of Rochester Plan to Cut Mathematics Is Recipe for
  Disaster}},
   publisher={The Chronicle of Higher Education},
        date={1996},
  url={https://www.chronicle.com/article/u-of-rochester-plan-to-cut-mathematics-is-recipe-for-disaster/},
  note={\url{https://www.chronicle.com/article/u-of-rochester-plan-to-cut-mathematics-is-recipe-for-disaster/},
  Accessed October 6, 2023},
}

\bib{calcscience}{article}{
      author={Kramer, Laird},
      author={Fuller, Edgar},
      author={Watson, Charity},
      author={Castillo, Adam},
      author={Oliva, Pablo~Duran},
      author={Potvin, Geoff},
       title={Establishing a new standard of care for calculus using trials
  with randomized student allocation},
        date={2023},
     journal={Science},
      volume={381},
      number={6661},
       pages={995\ndash 998},
      eprint={https://www.science.org/doi/pdf/10.1126/science.ade9803},
         url={https://www.science.org/doi/abs/10.1126/science.ade9803},
}

\bib{enrcliff2023}{misc}{
      author={Knox, Liam},
       title={Grasping for a foothold on the enrollment cliff},
   publisher={Inside Higher Ed},
        date={2023},
  url={https://www.insidehighered.com/news/business/revenue-strategies/2023/05/12/grasping-foothold-enrollment-cliff},
  note={\url{https://www.insidehighered.com/news/business/revenue-strategies/2023/05/12/grasping-foothold-enrollment-cliff},
  Accessed October 31, 2023},
}

\bib{Krebs2023}{misc}{
      author={Krebs, Paula},
       title={Opinion: {West Virginia University’s} cuts are a travesty},
   publisher={CNN},
        date={2023},
  url={https://edition.cnn.com/2023/09/19/opinions/west-virginia-university-gordon-gee-humanities-cuts-krebs/index.html},
  note={\url{https://edition.cnn.com/2023/09/19/opinions/west-virginia-university-gordon-gee-humanities-cuts-krebs/index.html},
  Accessed September 24, 2023},
}

\bib{axioms12070699}{article}{
      author={Kwessi, Eddy},
       title={Topological comparison of some dimension reduction methods using
  persistent homology on {EEG} data},
        date={2023},
        ISSN={2075-1680},
     journal={Axioms},
      volume={12},
      number={7},
         url={https://www.mdpi.com/2075-1680/12/7/699},
}

\bib{newR12023}{misc}{
      author={Lederman, Doug},
       title={{A New Approach to Categorizing Colleges}},
   publisher={Inside Higher Ed},
        date={2023},
  url={https://www.insidehighered.com/news/institutions/2023/11/01/major-overhaul-coming-key-framework-organizing-higher-ed},
  note={\url{https://www.insidehighered.com/news/institutions/2023/11/01/major-overhaul-coming-key-framework-organizing-higher-ed},
  Accessed November 1, 2023},
}

\bib{leyva2022making}{article}{
      author={Leyva, Elizabeth},
      author={Walkington, Candace},
      author={Perera, Harsha},
      author={Bernacki, Matthew},
       title={Making mathematics relevant: An examination of student interest
  in mathematics, interest in {STEM} careers, and perceived relevance},
        date={2022},
     journal={International Journal of Research in Undergraduate Mathematics
  Education},
      volume={8},
      number={3},
       pages={612\ndash 641},
}

\bib{malcom2020katherine}{article}{
      author={Malcom, Shirley~M},
       title={{Katherine Johnson (1918--2020)}},
        date={2020},
     journal={Science},
      volume={368},
      number={6491},
       pages={591\ndash 591},
}

\bib{pbs2019}{misc}{
      author={Marcus, Jon},
       title={Most {Americans} don’t realize state funding for higher ed fell
  by billions},
   publisher={Public Broadcasting System},
        date={2019},
  url={\url{https://www.pbs.org/newshour/education/most-americans-dont-realize-state-funding-for-higher-ed-fell-by-billions}},
  note={\url{https://www.pbs.org/newshour/education/most-americans-dont-realize-state-funding-for-higher-ed-fell-by-billions},
  Accessed September 8, 2023},
}

\bib{bobmoses2021}{misc}{
      author={Milloy, Courtland},
       title={{Bob Moses} saw math as the path to equality. {School} systems
  should build upon his work.},
   publisher={The Washington Post},
        date={2021},
  url={\url{https://www.washingtonpost.com/local/bob-moses-algebra-math-black-students/2021/07/27/74e41f24-eef5-11eb-81d2-ffae0f931b8f_story.html}},
        note={Published July 27, 2021 -
  \url{https://www.washingtonpost.com/local/bob-moses-algebra-math-black-students/2021/07/27/74e41f24-eef5-11eb-81d2-ffae0f931b8f_story.html},
  Accessed September 6, 2023},
}

\bib{unkept2018}{misc}{
      author={Mitchell, Michael},
      author={Leachman, Michael},
      author={Masterson, Kathleen},
      author={Waxman, Samantha},
       title={Unkept promises: State cuts to higher education threaten access
  and equity},
   publisher={Center on Budget and Policy Priorities},
        date={2018},
  url={https://www.cbpp.org/research/state-budget-and-tax/unkept-promises-state-cuts-to-higher-education-threaten-access-and},
  note={\url{https://www.cbpp.org/research/state-budget-and-tax/unkept-promises-state-cuts-to-higher-education-threaten-access-and},
  Accessed September 24, 2023},
}

\bib{ukrpk2022}{misc}{
      author={Moody, Josh},
       title={$16$ states underfunded land-grant {HBCUs} by over \$12b, {Biden}
  admin says},
   publisher={Inside Higher Ed},
        date={2023},
  url={https://www.insidehighered.com/news/2022/02/18/university-kansas-plans-cut-42-academic-programs},
  note={\url{https://www.insidehighered.com/news/2022/02/18/university-kansas-plans-cut-42-academic-programs},
  Accessed September 11, 2023},
}

\bib{menezes2018handbook}{book}{
      author={Menezes, Alfred~J},
      author={Van~Oorschot, Paul~C},
      author={Vanstone, Scott~A},
       title={Handbook of applied cryptography},
   publisher={CRC press},
        date={2018},
}

\bib{usnewsmath}{misc}{
      author={News, US},
      author={Report, World},
       title={Mathematician - overview},
        date={2023},
         url={\url{https://money.usnews.com/careers/best-jobs/mathematician}},
        note={\url{https://money.usnews.com/careers/best-jobs/mathematician},
  Accessed September 6, 2023},
}

\bib{sheff2023}{misc}{
      author={Officers, State Higher Education~Executive},
       title={State higher education finance: Fy 2022},
        date={2023},
         url={\url{https://shef.sheeo.org/data-downloads/}},
        note={\url{https://shef.sheeo.org/data-downloads/}, Accessed September
  6, 2023},
}

\bib{AceNetCC}{misc}{
      author={of~Institutions~of Higher~Education, The
  Carnegie~Classification},
       title={About {Carnegie} classification},
        date={2023},
        note={Retrieved September 4, 2023 from
  \url{https://carnegieclassifications.acenet.edu/}},
}

\bib{uncg2023}{misc}{
      author={of~North Carolina-Greensboro, University},
       title={Academic review dasboards},
   publisher={University of North Carolina-Greensboro},
        date={2018},
  url={https://innovation.uncg.edu/initiatives/academic-data-dashboards-admin-services-review/},
  note={\url{https://innovation.uncg.edu/initiatives/academic-data-dashboards-admin-services-review/},
  Accessed October 24, 2023},
}

\bib{olson2012engage}{article}{
      author={Olson, Steve},
      author={Riordan, Donna~Gerardi},
       title={Engage to excel: producing one million additional college
  graduates with degrees in science, technology, engineering, and mathematics.
  report to the president.},
        date={2012},
     journal={Executive Office of the President},
}

\bib{caks2022}{misc}{
      author={Perez, Suzanne},
       title={{Who needs college algebra? Kansas universities may rethink math
  requirements}},
   publisher={KMUW Wichita, KS},
        date={2022},
  url={https://www.kmuw.org/news/2022-12-12/who-needs-college-algebra-kansas-universities-may-rethink-math-requirements},
  note={\url{https://www.kmuw.org/news/2022-12-12/who-needs-college-algebra-kansas-universities-may-rethink-math-requirements},
  Accessed October 4, 2023},
}

\bib{quinn2023}{article}{
      author={Quinn, Ryan},
       title={{{WVU}, an {R-1} flagship…without language degrees or math
  {Ph.D.s}?}},
        date={2023},
     journal={Inside Higher Ed.},
        note={August 31, 2023, Retrieved from
  \url{https://www.insidehighered.com/news/institutions/research-universities/2023/08/31/wvu-r-1-flagshipwithout-language-degrees-or-math}},
}

\bib{rasmussen2019brief}{article}{
      author={Rasmussen, Chris},
      author={Apkarian, Naneh},
      author={Hagman, Jessica~Ellis},
      author={Johnson, Estrella},
      author={Larsen, Sean},
      author={Bressoud, David},
       title={Brief report: characteristics of precalculus through calculus 2
  programs: insights from a national census survey},
        date={2019},
     journal={Journal for Research in Mathematics Education},
      volume={50},
      number={1},
       pages={98\ndash 111},
}

\bib{reinholz2020time}{article}{
      author={Reinholz, Daniel~L},
      author={Rasmussen, Chris},
      author={Nardi, Elena},
       title={Time for (research on) change in mathematics departments},
        date={2020},
     journal={International Journal of Research in Undergraduate Mathematics
  Education},
      volume={6},
       pages={147\ndash 158},
}

\bib{mckuf}{misc}{
      author={Shanley, Garrett},
       title={{UF} signs \$4.7 million contract with global consulting firm},
   publisher={The Independent Florida Alligator},
        date={2023},
  url={\url{https://www.alligator.org/article/2023/08/uf-signs-4-7-million-contract-with-global-consulting-firm-mckinsey}},
        note={Published August 23, 2023 -
  \url{https://www.alligator.org/article/2023/08/uf-signs-4-7-million-contract-with-global-consulting-firm-mckinsey},
  Accessed September 6, 2023},
}

\bib{sassenyt2023}{misc}{
      author={Sokolove, Michael},
       title={{How Ben Sasse Became a Combatant in Florida’s Education
  Wars}},
   publisher={New York Times Magazine},
        date={2023},
  url={https://www.nytimes.com/2023/09/07/magazine/ben-sasse-university-florida.html},
  note={\url{https://www.nytimes.com/2023/09/07/magazine/ben-sasse-university-florida.html},
  Accessed September 11, 2023},
}

\bib{wvu2023}{misc}{
      author={Spitalniak, Laura},
       title={{WVU} board approves dramatic academic cuts to address {\$45m}
  deficit},
   publisher={Higher Ed Dive},
        date={2023},
  url={https://www.highereddive.com/news/wvu-board-approves-academic-cuts-45m-deficit/693857/},
  note={\url{https://www.highereddive.com/news/wvu-board-approves-academic-cuts-45m-deficit/693857/},
  Accessed September 19, 2023},
}

\bib{smith2017seminal}{inproceedings}{
      author={Smith, Wendy~M},
      author={Webb, David~C},
      author={Bowers, Janet},
      author={Voigt, Matt},
      author={Team, TS},
       title={{SEMINAL}: Preliminary findings on institutional changes in
  departments of mathematics},
organization={Association of Public and Land-grant Universities},
        date={2017},
   booktitle={Proceedings of the sixth annual mathematics teacher education
  partnership conference},
       pages={121\ndash 128},
}

\bib{chron2015}{misc}{
      author={Thomason, Andy},
       title={Another republican governor proposes big cuts in higher
  education},
   publisher={The Chronicle of Higher Education},
        date={2015},
  url={\url{https://www.chronicle.com/blogs/ticker/another-republican-governor-proposes-big-cuts-in-higher-education}},
        note={Published February 18, 2015 -
  \url{https://www.chronicle.com/blogs/ticker/another-republican-governor-proposes-big-cuts-in-higher-education},
  Accessed September 6, 2023},
}

\bib{wvuair2023}{misc}{
      author={Toney, Mike},
       title={{WVU expenses paid to aircraft charter provider have risen since
  Gee called for decreasing costs}},
   publisher={Charleston Gazette-Mail},
        date={2023},
  url={https://www.wvgazettemail.com/news/education/wvu-expenses-paid-to-aircraft-charter-provider-have-risen-since-gee-called-for-decreasing-costs/article_7a5d42eb-75b5-5216-b20b-0c8788cdc4fd.html},
  note={\url{https://www.wvgazettemail.com/news/education/wvu-expenses-paid-to-aircraft-charter-provider-have-risen-since-gee-called-for-decreasing-costs/article_7a5d42eb-75b5-5216-b20b-0c8788cdc4fd.html},
  Accessed October 6, 2023},
}

\bib{wvurec}{misc}{
      author={University, West~Virginia},
       title={Academic program portfolio review},
        date={2023},
        note={Retrieved September 13, 2023 from
  \url{https://provost.wvu.edu/academic-transformation/academic-program-portfolio-review?asd}},
}

\bib{fightwvu2023}{misc}{
      author={Warner, John},
       title={Fighting for the future at {WVU}},
   publisher={Inside Higher Ed},
        date={2023},
  url={https://www.insidehighered.com/opinion/blogs/just-visiting/2023/09/15/change-coming-higher-ed-who-will-lead-it},
  note={\url{https://www.insidehighered.com/opinion/blogs/just-visiting/2023/09/15/change-coming-higher-ed-who-will-lead-it},
  Accessed September 24, 2023},
}

\bib{nytwvu}{misc}{
      author={Weatherby, Leif},
       title={What just happened at {West Virginia University} should worry all
  of us},
   publisher={The New York Times},
        date={2023},
  url={\url{https://www.nytimes.com/2023/08/20/opinion/west-virginia-university-cuts.html}},
        note={Published August 20, 2023 -
  \url{https://www.nytimes.com/2023/08/20/opinion/west-virginia-university-cuts.html},
  Accessed September 6, 2023},
}

\end{biblist}
\end{bibdiv}

\end{document}